\documentclass[12pt]{article} 
\usepackage{amsfonts,latexsym,color}
\input{psfig.sty}

\newcounter{environment}[section]
\renewcommand{\theenvironment}{%
\arabic{section}.\arabic{environment}}

{\begin{rm}\refstepcounter{environment}{\textbf\theenvironment\
\bf Definition.~~}}%
{\end{rm}}

{\begin{rm}\refstepcounter{environment}{\textbf\theenvironment\
\bf Example.~~}}%
{\end{rm}}

{\begin{rm}\refstepcounter{environment}{\textbf\theenvironment\
\bf Conjecture.~~}}%
{\end{rm}}

{\begin{rm}\refstepcounter{environment}{\textbf\theenvironment\
\bf Proposition.~~}}%
{\end{rm}}

{\begin{rm}\refstepcounter{environment}{\textbf\theenvironment\
\bf Proposition and Definition.~~}}%
{\end{rm}}

\newenvironment{theorem}%
{\begin{rm}\refstepcounter{environment}{\textbf\theenvironment\
\bf Theorem.~~}}%
{\end{rm}}

\newenvironment{corollary}%
{\begin{rm}\refstepcounter{environment}{\textbf\theenvironment\
\bf Corollary.~~}}%
{\end{rm}}

{\begin{rm}\refstepcounter{environment}{\textbf\theenvironment\
\bf Lemma.~~}}%
{\end{rm}}

\begin{document}
\newcommand{\qbc}[2]{{\left [{#1 \atop #2} \right ]}}
\newcommand{\df}[2]{1-x_{#1}x_{#2}^{-1}}
\newcommand{\pdf}[2]{\left(1-x_{#1}x_{#2}^{-1}\right)}
\newcommand{\beq}{\begin{equation}}
\newcommand{\eeq}{\end{equation}}
\newcommand{\be}{\begin{enumerate}}
\newcommand{\ee}{\end{enumerate}}
\newcommand{\bea}{\begin{eqnarray}}
\newcommand{\eea}{\end{eqnarray}}
\newcommand{\beas}{\begin{eqnarray*}}
\newcommand{\eeas}{\end{eqnarray*}}
\newcommand{\ds}{\displaystyle}
\newcommand{\rr}{\mathbb{R}}
\newcommand{\zz}{\mathbb{Z}}
\newcommand{\nn}{\mathbb{N}}
\newcommand{\pp}{\mathbb{PS}}
\newcommand{\st}{\,:\,}
\newcommand{\ps}{\mathbb{Q}[[x]]_0}
\newcommand{\ii}{\mathrm{in}}
\newcommand{\ff}{\mathrm{fin}}
\newcommand{\cp}{{\cal P}}
\newcommand{\cf}{{\cal F}}
\newcommand{\sn}{\mathfrak{S}_n}
\newcommand{\fs}{\mathfrak{S}}
\newcommand{\vk}{\stackrel{\mathrm{TVK}}{\longrightarrow}}


\begin{centering}
\textcolor{red}{\Large\bf On the Enumeration}\\
\textcolor{red}{\Large\bf of Skew Young Tableaux}\\[.2in] 
\textcolor{blue}{Richard P. Stanley}\footnote{Partially supported by
  NSF grant 
\#DMS-9988459.}\\ 
Department of Mathematics\\
Massachusetts Institute of Technology\\
Cambridge, MA 02139\\
\emph{e-mail:} rstan@math.mit.edu\\[.2in]
\textcolor{magenta}{version of 9 September 2001}\\[.2in]

\end{centering}
\vskip 10pt
\section{Introduction.} A recent paper \cite{m-m-w} of McKay, Morse,
and Wilf considers the number $N(n;T)$ of standard Young tableaux
(SYT) with $n$ cells that contain a fixed standard Young tableau $T$
of shape $\alpha\vdash k$. (For notation and terminology related to
symmetric functions and tableaux, see \cite{macd} or \cite[Ch.\
7]{ec2}.) They obtain the asymptotic formula
  \beq N(n;T) \sim \frac{t_nf^\alpha}{k!}, \label{eq:mmw} \eeq
where $f^\alpha$ denotes the number of SYT of shape $\alpha$ and
$t_n$ denotes the number of involutions in the symmetric group $\sn$. 
Note that $N(n;T)=N(n;U)$ whenever $T$ and $U$ are SYT of the same
shape. Hence we can write $N(n;\alpha)$ for $N(n;T)$. Moreover, it is
clear that
  \beq N(n;\alpha) = \sum_{\lambda\vdash n}f^{\lambda/\alpha},
  \label{eq:nnt} \eeq
where $f^{\lambda/\alpha}$ denotes the number of SYT $T$ of skew
shape $\lambda/\alpha$.

In Section \ref{sec2} we extend equation (\ref{eq:mmw}), using
techniques from the theory of symmetric functions, to give an explicit
formula for $N(n;\alpha)$ as a finite linear combination of 
$t_{n-j}$'s, from which in principle we can write down the entire
asymptotic expansion of $N(n;\alpha)$. In Section \ref{sec3} we apply
similar techniques, together with asymptotic formulas for character
values of $\sn$ due to Biane \cite{biane} and to Vershik and Kerov
\cite{v-k}, to derive the asymptotic behavior of $f^{\lambda/\alpha}$
as a function of $\lambda$ for fixed $\alpha$. 

\section{A formula for $N(n;\alpha)$.} \label{sec2}
Let $\chi^\alpha(\lambda)$ denote the value of the irreducible
character $\chi^\alpha$ of $\fs_k$ on a permutation of cycle type
$\lambda\vdash k$ (as explained e.g.\ in \cite[{\S}1.7]{macd} or
\cite[{\S}{\S}7.17--7.18]{ec2}). Let $m_i(\mu)$ denote the number of
parts of the partition $\mu$ equal to $i$, and write $\tilde{\mu}$ for
the partition obtained from $\mu$ by replacing every even part $2i$
with the two parts $i,i$. For instance,
  $$ \mu=(6,6,5,4,2,1)\ \Rightarrow\ \tilde{\mu}=
       (5,3,3,3,3,2,2,1,1,1). $$
Equivalently, if $w$ is a permutation of cycle type $\mu$, then $w^2$
has cycle type $\tilde{\mu}$. Note that a permutation of cycle type
$\tilde{\mu}$ is necessarily even. We will use notation such as
$(\tilde{\mu},1^{k-j})$ to denote a partition whose parts are the
parts of $\tilde{\mu}$ with $k-j$ additional parts equal to 1. Finally
we let $z_\mu$ denote the number of permutations commuting with a
fixed permutation of cycle type $\mu$, so 
  $$ z_\mu = 1^{m_1(\mu)}2^{m_2(\mu)}\cdots m_1(\mu)!\,
        m_2(\mu)!\cdots. $$

The main result of this section is the following.

\begin{theorem} \label{thm1}
\emph{Let $\alpha\vdash k$. Then for $n\geq k$ we have}
  \beq N(n;\alpha) = \sum_{j=0}^k \frac{t_{n-j}}{(k-j)!}
  \sum_{{\mu\vdash j\atop  m_1(\mu)=m_2(\mu)=0}}
  z_\mu^{-1}\chi^\alpha(\tilde{\mu},1^{k-j}). \label{eq:thm1} \eeq
 \end{theorem}

\textbf{Proof.} Let $\lambda\vdash n\geq k$, and let
$p_\lambda=p_{\lambda_1}p_{\lambda_2}\cdots$ denote the power sum
symmetric function indexed by $\lambda$. Similarly
$s_{\lambda/\alpha}$ denotes the skew Schur function indexed by
$\lambda/\alpha$. Since for any homogeneous symmetric function $f$ of
degree $n-k$ we have that $\langle p_1^{n-k},f\rangle$ is the
coefficient of $x_1\cdots x_{n-k}$ in $f$, and since the coefficient
of $x_1\cdots x_{n-k}$ in $s_{\lambda/\alpha}$ is $f^{\lambda/\alpha}$,
we have (using a basic property \cite[(5.1)]{macd}\cite[Thm.\
7.15.4]{ec2} of the standard scalar product 
$\langle \cdot,\cdot\rangle$ on symmetric functions)
  \beas f^{\lambda/\alpha} & = & \langle
   p_1^{n-k},s_{\lambda/\alpha}\rangle \\ & = &
   \langle p_1^{n-k}s_\alpha,s_\lambda\rangle. \eeas
Summing on $\lambda\vdash n$ gives
  \beq N(n;\alpha)=\left\langle p_1^{n-k}s_\alpha, \sum_{\lambda\vdash
      n} s_\lambda\right\rangle. \label{eq:nna} \eeq
Now \cite[Exam.\ I.5.4, p.\ 76]{macd}\cite[Cor.\ 7.13.8]{ec2}
  $$ \sum_\lambda s_\lambda = \frac{1}{\prod_i (1-x_i)\cdot
     \prod_{i<j}(1-x_ix_j)}, $$
summed over all partitions $\lambda$ of all $n\geq 0$. Since
  \beas \frac{1}{\prod_i (1-x_i)\cdot
     \prod_{i<j}(1-x_ix_j)} & = & \exp\sum_{n\geq 1}\frac 1n
       \left(\sum_i x_i^n + \sum_{i<j}x_i^nx_j^n\right)\\
      & = & \exp \left(\sum_{n\geq 1}\frac{p_{2n-1}}{2n-1}
           +\sum_{n\geq 1}\frac{p_n^2}{2n}\right), \eeas
there follows
  \bea \sum_{\lambda\vdash n}s_\lambda & = & \sum_{\lambda=(1^{m_1},
      2^{m_2},\dots)\vdash n} z_\lambda^{-1}p_1^{m_1+2m_2}
       p_2^{2m_4}p_3^{m_3+2m_6}p_4^{2m_8}\cdots\nonumber\\
     & = & \sum_{\lambda\vdash n}z_\lambda^{-1}p_{\tilde{\lambda}},
      \label{eq:ssum} \eea
where $(1^{m_1},2^{m_2},\dots)$ denotes the partition with $m_i$
parts equal to $i$. 

It follows from \cite[p.\ 76]{macd}\cite[solution to Exer.\
7.35(a)]{ec2} that for any symmetric functions $f$ and $g$ we have 
  $$ \langle p_1 f,g\rangle = \left\langle f, 
   \frac{\partial}{\partial p_1}g\right\rangle, $$
where $\frac{\partial}{\partial p_1}g$ indicates that we are to expand
$g$ as a polynomial in the $p_i$'s and then differentiate with respect
to $p_1$. Applying this to equation~(\ref{eq:nna}) and using
(\ref{eq:ssum}) yields
  \bea N(n;\alpha) & = & \left\langle s_\alpha, \frac{\partial^{n-k}}
  {\partial p_1^{n-k}} \sum_{\lambda
      \vdash n} z_\lambda^{-1}p_{\tilde{\lambda}}\right\rangle
    \nonumber\\
    & = &  \left\langle s_\alpha,\sum_{\lambda\vdash n} z_\lambda^{-1}
       (m_1+2m_2)_{n-k}\,p_1^{-n+k}p_{\tilde{\lambda}}\right\rangle, 
     \label{eq:m1m2} \eea
where $m_i=m_i(\lambda)$ and $(a)_{n-k}=a(a-1)\cdots (a-n+k+1)$. 

Fix $m_1+2m_2=n-j$ in equation (\ref{eq:m1m2}). Thus $\lambda =
(\mu,2^{m_2},1^{m_1})$ for some unique $\mu\vdash k$ satisfying
$m_1(\mu)= 
m_2(\mu)=0$. Since $n!/z_\lambda$ is the number of permutations in
$\sn$ of cycle type $\lambda$, we have for fixed $\mu\vdash j$ that
  $$ \sum_{\lambda=(\mu,2^{m_2},1^{m_1})} \frac{n!}{z_\lambda} = 
      t_{n-j}{n\choose j}\frac{j!}{z_\mu}. $$
Moreover, 
  $$ p_1^{-n+k}p_{\tilde{\lambda}}=p_1^{k-j}p_{\tilde{\mu}}. $$
It follows that 
  \beas N(n;\alpha) & = & \left\langle s_\alpha, \sum_{j=0}^k     
      j!{n\choose j}(n-j)_{n-k}\frac{t_{n-j}}{n!} 
     \sum_{{\mu\vdash j\atop m_1(\mu)=m_2(\mu)=0}} z_\mu^{-1}
       p_1^{k-j}p_{\tilde{\mu}} \right\rangle\\ & = &
       \left\langle s_\alpha, \sum_{j=0}^k     
      \frac{t_{n-j}}{(k-j)!} 
     \sum_{{\mu\vdash j\atop m_1(\mu)=m_2(\mu)=0}} z_\mu^{-1}
       p_1^{k-j}p_{\tilde{\mu}} \right\rangle. \eeas
Since \cite[(7.7)]{macd}\cite[p.\ 348]{ec2}
  $$ \langle s_\alpha,p_1^{k-j}p_{\tilde{\mu}}\rangle =
       \chi^{\alpha}(\tilde{\mu},1^{k-j}), $$
the proof follows. $\ \Box$

Note that the restriction $n\geq k$ in Theorem~\ref{thm1} is
insignificant since $N(n;\alpha)=0$ for $n<k$.

Theorem~\ref{thm1} expresses $N(n;\alpha)$ as a linear
combination of the functions $t_{n-j}$, $0\leq j\leq k$. Since
$t_{n-j-1} =o(t_{n-j})$, this formula for $N(n;\alpha)$ is actually
an asymptotic expansion. The first few terms are
  \bea N(n;\alpha) & = & \frac{1}{k!}f^\alpha t_n + \frac{1}{3
    (k-3)!} \chi^\alpha(3,1^{k-3})t_{n-3}\nonumber\\ & &
  \ + \frac{1}{4 (k-4)!}
      \chi^\alpha(2,2,1^{k-4})t_{n-4}+\frac{1}{5(k-5)!}
      \chi^\alpha(5,1^{k-5})t_{n-5}\nonumber\\ & & \ + \frac{2}{9
        (k-6)!} \chi^{\alpha}(3,3,1^{k-6})t_{n-6} + O(t_{n-7}). 
     \label{eq:nnaas} \eea
Note that by symmetry it is clear that if $\alpha'$ is the conjugate
partition to $\alpha$ then $N(n;\alpha)=N(n;\alpha')$. Indeed, since a
permutation of cycle type $\tilde{\mu}$ is even we have
$\chi^\alpha(\tilde{\mu},1^{k-j})
=\chi^{\alpha'}(\tilde{\mu},1^{k-j})$. The exact formulas for 
$N(n;\alpha)$ when $|\alpha|\leq 5$ and $|\alpha|\leq n$ are given as
follows (where we write e.g.\ $N(n;21)$ for $N(n;(2,1))$):
  \beas N(n;1) & = & t_n\\
       N(n;2)\ =\ N(n;11) & = & \frac 12 t_n\\
       N(n;3)\ =\ N(n;111) & = & \frac 16(t_n+2t_{n-3})\\
       N(n;21) & = & \frac 13(t_n-t_{n-3})\\
     N(n;4)\ =\ N(n;1111) & = & \frac{1}{24}(t_n+8t_{n-3}+6t_{n-4})\\
     N(n;31)\ =\ N(n;211) & = & \frac 18(t_n-2t_{n-4})\\
     N(n;22) & = & \frac{1}{12}(t_n-4t_{n-3}+6t_{n-4})\\
     N(n;5)\ =\ N(n;11111) & = & \frac{1}{120}(t_n+20t_{n-3}+30t_{n-4}
                  +24t_{n-5})\\
     N(n;41)\ =\ N(n;2111) & = & \frac{1}{30}(t_n+5t_{n-3}-6t_{n-5})\\
     N(n;32)\ =\ N(n;221) & = & \frac{1}{24}(t_n-4t_{n-3}+6t_{n-4})\\
     N(n;311) & = & \frac{1}{20}(t_n-10t_{n-4}+4t_{n-5}).    
  \eeas
\indent The complete asymptotic expansion of $t_n$ beginning
  $$ t_n \approx \frac{1}{\sqrt{2}}n^{n/2}e^{-\frac n2+\sqrt{n}-\frac 14}
      \left( 1+\frac{7}{24\sqrt{n}}-
    \frac{119}{1152n}+\cdots\right) $$
was obtained by Moser and Wyman \cite[3.39]{m-w}. In principle this
can be used to obtain the asymptotic expansion of $N(n;\alpha)$ in
terms of more ``familiar'' functions than $t_{n-j}$. The first few
terms can be obtained from the formula
  $$ t_{n-j} = \frac{1}{\sqrt{2}}n^{\frac{n-j}{2}}
        e^{-\frac n2+\sqrt{n}-\frac 14}\left( 1 +\left( \frac{7}{24}
         -\frac
         j2\right)\frac{1}{\sqrt{n}}
     -\left( \frac{119}{1152}+\frac{7}{48}j-\frac 38 j^2\right)
       \frac 1n \right. $$
  $$ \qquad\qquad\left. +O\left( \frac{1}{n^{3/2}}\right)\right), $$
though we omit the details.

Instead of counting the number $N(n;\alpha)$ of SYT with $n$
cells containing a fixed SYT $T$ of shape $\alpha$, we can ask (as
also done in \cite{m-m-w}) for the probability $P(n;\alpha)$ that a
random SYT with $n$ cells (chosen from the uniform distribution on all
SYT with $n$ cells) contains $T$ as a subtableau. Since the total
number of SYT with $n$ cells is $t_n$, we have 
  $$ P(n;\alpha) = N(n;\alpha)/t_n. $$
Let $e_j(\alpha)$ denote the coefficient of $t_{n-j}$ in the
right-hand side of (\ref{eq:thm1}), viz.,
  \beq e_j(\alpha) = \frac{1}{(k-j)!}
  \sum_{{\mu\vdash j\atop  m_1(\mu)=m_2(\mu)=0}}
  z_\mu^{-1}\chi^\alpha(\tilde{\mu},1^{k-j}). \label{eq:ej} \eeq
It follows from Theorem~\ref{thm1}, using the fact that $e_0(\alpha)
=f^\alpha/k!$ and $e_1(\alpha)=e_2(\alpha)=0$, that
  $$ P(n;\alpha) = \frac{f^\alpha}{k!}+\frac{e_3(\alpha)}{n^{3/2}}
         -\frac{3e_3(\alpha)-2e_4(\alpha)}{n^2}+O\left(
           n^{-5/2}\right). $$ 
The leading term of this expansion was obtained in \cite[Thm.\
1]{m-m-w}. 

There is an alternative formula for $N(n;\alpha)$ which, though not as
convenient for asymptotics, is more combinatorial than equation
(\ref{eq:thm1}) because it avoids using the characters of $\sn$. This
formula could be derived directly from Theorem~\ref{thm1}, but we give
an alternative proof which is implicitly bijective (since the formulas
on which it is based have bijective proofs).

\begin{theorem} \label{thmnew}
\emph{Let $\alpha\vdash k$. Then for all $n\geq 0$ we have}
  \beq N(n+k;\alpha)= \sum_{j=0}^k {n\choose j}\left( \sum_{\mu\vdash k-j}
    f^{\alpha/\mu}\right) t_{n-j}. \label{eq:new} \eeq
\end{theorem}

\textbf{Proof.} We begin with the following Schur function identity,
proved independently by Lascoux, Macdonald, Towber, Stanley,
Zelevinsky, and perhaps others. This identity appears in \cite[Exam.\ 
I.5.27(a), p.\ 93]{macd}\cite[Exer.\ 7.27(e)]{ec2} and was given a
bijective proof by Sagan and Stanley \cite[Cor.\ 6.4]{s-s}:
  $$ \sum_{\lambda} s_{\lambda/\alpha} = \frac{1}
  {\prod_i(1-x_i)\cdot \prod_{i<j}(1-x_ix_j)}\sum_\mu
  s_{\alpha/\mu}. $$
Apply the homomorphism ex that takes the power sum symmetric function
$p_n$ to $\delta_{1n}u$, where $u$ is an indeterminate. This
homomorphism is the exponential specialization discussed in \cite[pp.\ 
304--305]{ec2}. Two basic properties of ex are the following:
  $$ \mathrm{ex}(f) = \sum_{n\geq 0}[x_1x_2\cdots x_n]f
      \frac{u^n}{n!} $$ 
  $$ \mathrm{ex}\frac{1} {\prod_i(1-x_i)\cdot \prod_{i<j}(1-x_ix_j)} = 
      e^{u+\frac 12u^2}, $$
where $[x_1x_2\cdots x_n]f$ denotes the coefficient of $x_1x_2\cdots
x_n$ in $f$. Since 
  $$ [x_1x_2\cdots x_n]s_{\lambda/\alpha}=
  f^{\lambda/\alpha},\ \mathrm{when}\ |\lambda/\alpha|=n, $$ 
we obtain
  \beq \sum_{n\geq 0} \frac{u^n}{n!}\sum_{\lambda\vdash n+k}
     f^{\lambda/\alpha} = e^{u+\frac 12u^2}\sum_{j=0}^k
    \frac{u^j}{j!}\sum_{\mu\vdash k-j}f^{\alpha/\mu}. 
   \label{eq:fcor} \eeq
Taking the coefficient of $u^n/n!$ on both sides yields
(\ref{eq:new}). $\ \Box$

\begin{corollary} \label{corgf}
\emph{We have}
  $$ \sum_{n\geq 0} \sum_\alpha N(n+|\alpha|;\alpha)s_\alpha
    \frac{u^n}{n!} = \left( \sum_\mu s_\mu\right) e^{(p_1+1)u+\frac
    12 u^2}. $$
\end{corollary}

\textbf{Proof.} Multiply (\ref{eq:fcor}) by $s_\alpha$ and sum on
$\alpha$ to get
  \beas \sum_{n\geq 0} \sum_\alpha N(n+|\alpha|;\alpha)s_\alpha
    \frac{u^n}{n!} & = & e^{u+\frac 12 u^2}\sum_{j\geq 0} 
       \frac{u^j}{j!}
        \sum_{|\alpha/\mu|=j} f^{\alpha/\mu}s_\alpha\\ & = &
   e^{u+\frac 12 u^2}\sum_{j\geq 0} \frac{u^j}{j!}
        \sum_{|\alpha/\mu|=j}\langle p_1^j,s_{\alpha/\mu}
    \rangle s_\alpha\\ & = &
   e^{u+\frac 12 u^2}\sum_{j\geq 0} \frac{u^j}{j!} 
        \sum_{|\alpha/\mu|=j}\langle p_1^j s_\mu,s_\alpha
    \rangle s_\alpha\\ & = & 
    e^{u+\frac 12 u^2}\sum_{j\geq 0} \frac{u^j}{j!} \sum_\mu 
       p_1^j s_\mu\\ & = & 
     \left( \sum_\mu s_\mu\right) e^{(p_1+1)u+\frac 12 u^2}.
    \ \Box \eeas

The case when $\alpha$ consists of a single row (or column) is
particularly simple, since then each
$\chi^\alpha(\tilde{\mu},1^{k-j})=1$ in (\ref{eq:thm1}). We will then
write $N(n;k)$ as short for $N(n;(k))$. The
coefficient $e_j(\alpha)$ becomes simply
$e_j(k)=q_j/(k-j)!$, where $j!q_j$ is the number of permutations
$w\in\sn$ with no cycles of length one or two. By standard enumerative
reasoning (see e.g.\ \cite[Exam.\ 5.2.10]{ec2}) we have
  \beq \sum_{j\geq 0} q_j x^j = \frac{e^{-x-\frac 12
      x^2}}{1-x}. \label{eq:qj} \eeq
From this and Theorems~\ref{thm1} and \ref{thmnew} it is easy to
deduce the following results, which we simply state without proof.

\begin{corollary} \label{cor:abc}
(a) \emph{We have} 
  $$ N(n+k;k) = \sum_{j=0}^k {n\choose j}t_{n-j}
     = \sum_{j=0}^k\frac{q_j}{(k-j)!}t_{n+k-j}, $$
\emph{where $q_j$ is given by (\ref{eq:qj}).}

(b) \emph{Define polynomials $A_n(x)$ by $A_0(x)=1$ and}
 $$ A_{n+1}(x) = A'_n(x)+(x+1)A_n(x),\quad n\geq 0. $$
\emph{Then}
  $$ \sum_{k\geq 0}  N(n+k;k)x^k = \frac{A_n(x)}{1-x}. $$
(c) \emph{Let} 
  $$ e^{\frac 12u^2+2u}=\sum_{n\geq 0}b_n\frac{u^n}{n!}. $$
\emph{Then $N(n+k;k)=b_n$ if $n\leq k$.} 
\end{corollary}
 
The stability property of Corollary~\ref{cor:abc}(c) is easy to see
by direct combinatorial reasoning. If $n\leq k$, then a skew SYT of
shape $\lambda/\alpha$, where $\lambda\vdash n+k$ and $\alpha\vdash
k$,  consists of a first row containing some $j$-element subset of
$1,2,\dots,n$, together with some disjoint SYT $U$ on the remaining
$n-j$ letters. There are $t_{n-j}$ possibilities for $U$, so
  $$ N(n+k;k) = \sum_{j=0}^n {n\choose j}t_{n-j}, $$
which is equivalent to Corollary~\ref{cor:abc}(c).
  
\section{Asymptotics of $f^{\lambda/\alpha}$.} \label{sec3} 
Rather than considering the sum $\sum_{\lambda\vdash n}
f^{\lambda/\alpha}$, we could investigate instead the individual
terms $f^{\lambda/\alpha}$. The analogue of Theorem~\ref{thm1} is the
following. 

\begin{theorem} \label{thm2}
\emph{Let $\alpha\vdash k$ and $n\geq k$. Then for any partition
  $\lambda\vdash n$ we have}  
  \beq f^{\lambda/\alpha} = \sum_{\nu\vdash k} z_\nu^{-1}
       \chi^\lambda(\nu,1^{n-k})\chi^\alpha(\nu). \label{eq:fla} \eeq
\end{theorem}

\textbf{Proof.} The proof parallels that of
Theorem~\ref{thm1}. Instead of the power sum expansion of
$\sum_{\lambda\vdash n} s_\lambda$, we need the expansion of
$s_\lambda$ (where $\lambda\vdash n$), given by
\cite[p.\ 114]{macd}\cite[Cor.\ 7.17.5]{ec2} 
  $$ s_\lambda =\sum_{\mu\vdash n}z_\mu^{-1}\chi^\lambda(\mu)p_\mu. $$   
We therefore have
  \beas f^{\lambda/\alpha} & = & \left\langle p_1^{n-k},
     s_{\lambda/\alpha}\right\rangle\\ & = &
    \left\langle s_\alpha p_1^{n-k}, \sum_{\mu\vdash k}z_\mu^{-1}
    \chi^\lambda(\mu)p_\mu\right\rangle\\ & = &
     \left\langle s_\alpha,\frac{\partial^{n-k}}{\partial p_1^{n-k}}
     \sum_{\mu\vdash n}z_\mu^{-1}
    \chi^\lambda(\mu)p_\mu\right\rangle\\ & = &
   \left\langle s_\alpha,
     \sum_{\nu\vdash k}z_{(\nu,1^{n-k})}^{-1}\chi^\lambda(\nu,1^{n-k})
     (n-k+m_1(\nu))_n\, p_\nu\right\rangle\\ & = &
    \sum_{\nu\vdash k}z_{(\nu,1^{n-k})}^{-1}(n-k+m_1(\nu))_n\,  
    \chi^\lambda(\nu,1^{n-k})\chi^\alpha(\nu). \eeas
But
   $$ z_{(\nu,1^{n-k})}^{-1}(n-k+m_1(\nu))_n = z_\nu^{-1}, $$
and the proof follows. $\ \Box$

Theorem~\ref{thm2} can also be proved by inverting the formula given
in \cite[Exer.\ 7.62]{ec2}.

We would like to regard equation (\ref{eq:fla}) as an asymptotic
formula for $f^{\lambda/\alpha}$ when $\alpha$ is fixed and $\lambda$
is ``large.'' For this we need an asymptotic formula for
$\chi^\lambda(\nu,1^{n-k})$ when $\nu$ is fixed. Such a formula will
depend on the way in which the partitions $\lambda$ increase. The
first condition considered here is the following. Let
$\lambda^1,\lambda^2, \dots$ be a sequence of partitions such that
$\lambda^n\vdash n$, and such that the diagrams of the $\lambda^n$'s,
rescaled by a factor $n^{-1/2}$ (so that they all have area one)
converge uniformly to some limit $\omega$. (See \cite{biane} for a
more precise statement.) We will denote this convergence by
$\lambda^n\rightarrow \omega$. The following result is due to Biane
\cite{biane}, building on work of Vershik and Kerov.

\begin{theorem} \label{thm:biane}
\emph{Suppose that $\lambda^n\rightarrow \omega$. Then for $i\geq 2$
there exist constants (defined explicitly in \cite{biane}) 
$C_i(\omega)$, with $C_2(\omega)=1$, such that for any fixed partition
$\nu\vdash k$ of length $\ell(\nu)$ we have}
  $$ \chi^{\lambda^n}(\nu,1^{n-k}) =
  f^{\lambda^n}\left(\prod_{i=0}^{\ell(\nu)} C_{\nu_i+1}(\omega)\right)
  n^{-\frac 12(k-\ell(\nu))}\left( 1+O(1/n)\right), $$
\emph{as $n\rightarrow\infty$.}
\end{theorem}

Let $c_\nu=z_\nu^{-1}\chi^\lambda(\nu,1^{n-k})\chi^\alpha(\nu)$. 
It follows from Theorem~\ref{thm:biane} that $c_{(21^{k-2})}=
O\left( c_{(1^k)}n^{-1/2}\right)$, while $c_\nu=O\left( c_{(1^k)}n^{-1}
\right)$ and $c_\nu=O\left( c_{(21^{k-2})}n^{-1/2}\right)$ for
$\ell(\nu) \leq k-2$. Hence if $\lambda^n\rightarrow \omega$ then 
  \bea f^{\lambda^n/\alpha} & = &
  \left(z_{(1^k)}^{-1}\chi^{\lambda^n}(1^n)\chi^\alpha(1^k)
  +z_{(21^{k-2})}^{-1}\chi^{\lambda^n}(21^{n-2})\chi^\alpha(21^{k-2})
 \right)(1+O(1/n))\nonumber \\ & = &
   f^{\lambda^n}\left( \frac{1}{k!}f^\alpha+\frac{1}{2(k-2)!} 
   C_3(\omega)\chi^\alpha(21^{k-2})\frac{1}{\sqrt{n}}
     +O(1/n)\right). \label{eq:bias}\eea
Let us note that by \cite[p.\ 118]{macd}\cite[Exer.\ 7.51]{ec2} the
integer $\chi^\alpha(21^{k-2})$ appearing in (\ref{eq:bias}) has the
explicit value
  $$ \chi^\alpha(21^{k-2})=f^\alpha\frac{\sum {\alpha_i\choose 2}
   -\sum {\alpha'_i\choose 2}}{{k\choose 2}}. $$

The leading term of the right-hand side of (\ref{eq:bias}) is
independent of $\omega$, and in fact it follows 
from \cite{biane} that $f^{\lambda^n/\alpha}\sim \frac{1}{k!} f^\alpha
f^{\lambda^n}$ holds under the weaker
hypothesis that there exists a constant $A>0$ for which
$\lambda_1^n< A\sqrt{n}$ and $\ell(\lambda^n)<A\sqrt{n}$ for all
$n\geq 1$. 

\pagebreak
Given $\epsilon>0$, let 
  $$ \mathrm{Par}_\epsilon(n)=\{ \lambda\vdash n\st 
    (2-\epsilon)\sqrt{n}<\lambda_1<(2+\epsilon)\sqrt{n} $$ 
\vspace{-.3in}
  $$ \qquad \mathrm{and}\ (2-\epsilon)\sqrt{n}<\ell(\lambda)<
  (2+\epsilon)\sqrt{n}\}. $$  
It is a consequence of the work of Logan and Shepp \cite{l-s} or
Vershik and Kerov \cite{v-k2} (see e.g.\ \cite{b-d-j} for much stronger
results) that for any $\epsilon>0$, 
  $$ \sum_{\lambda\in\mathrm{Par}_\epsilon(n)} f^\lambda \sim t_n,\ \
     n\rightarrow\infty. $$
Thus not only is the sum $N(n;\alpha) =
\sum_{\lambda\vdash n}f^{\lambda/\alpha}$ asymptotic to $f^\alpha
t_n/k!$ as $n\rightarrow \infty$ (as follows from (\ref{eq:nnaas})),
but the terms $f^{\lambda/\alpha}$ contributing to ``most'' of the sum
are ``close'' to $f^\alpha f^\lambda/k!$.  

Another way of letting $\lambda$ become large was
considered by Vershik and Kerov in \cite{v-k} and in many subsequent
papers (after first being introduced by Thoma). Let
$\lambda^1,\lambda^2,\dots$ be a sequence of partitions 
such that $\lambda^n\vdash n$ and such that for all $i>0$, there exist
real numbers $a_i\geq 0$ and 
$b_i\geq 0$ satisfying $\sum_i(a_i+b_i)=1$ and
  \beas \lim_{n\rightarrow\infty}\frac{\lambda_i^n}{n} & = & a_i\\[.05in]
    \lim_{n\rightarrow\infty}\frac{\left(\lambda^n\right)'_i}{n} & = &
    b_i, \eeas
where $\left(\lambda^n\right)'_i$ denotes the $i$th part of the
conjugate partition to $\lambda^n$ (i.e., the length of the $i$th
column of the diagram of $\lambda^n$). We denote this situation by
$\lambda^n \vk (a;b)$, where
$a=(a_1,a_2,\dots)$ and $b=(b_1,b_2,\dots)$. For instance, if
$\lambda^{2n}= (n,n)$ and $\lambda^{2n-1}=(n,n-1)$, then $\lambda^n\vk
((1/2,1/2,0,\dots);(0,0,\dots))$. 

The following result is immediate from \cite{v-k}.

\begin{theorem} \label{thm:v-k}
\emph{Let $\lambda^n\vk (a;b)$. Then for any fixed partition
  $\nu\vdash k$,}   
  $$ \chi^{\lambda^n}(\nu, 1^{n-k}) =
        f^{\lambda^n}\prod_{j=1}^{\ell(\nu)}\left( \sum_i
          \alpha_i^{\nu_j} + (-1)^{\nu_j-1}\sum_i
          \beta_i^{\nu_j}\right) \left( 1+O(1/n)\right). $$
\end{theorem}

It follows that from Theorems~\ref{thm2} and \ref{thm:v-k} we have for
fixed $\alpha\vdash k$ the asymptotic formula
  \beq f^{\lambda^n/\alpha} = f^{\lambda^n} \left[\sum_{\nu\vdash k}
    z_\nu^{-1} \prod_{j=1}^{\ell(\nu)}\left( \sum_i
          \alpha_i^{\nu_j} + (-1)^{\nu_j-1}\sum_i
          \beta_i^{\nu_j}\right)\right](1+O(1/n)). \label{eq:flasum} \eeq

Now let $s_\lambda(x\,/\,y)$ denote the super-Schur function
indexed by $\lambda\vdash n$ in the variables
$x=(x_1,x_2,\dots)$ and $y=(y_1,y_2,\dots)$ \cite[Exam.\ 
I.23--I.24]{macd}, defined by $s_\lambda(x\,/\,-y) =\omega_y
s_\lambda(x,y)$ (where $\omega_y$ denotes the standard
involution $\omega$ acting on the $y$-variables only). (Note
that our $s_\lambda(x\,/\,y)$ corresponds to
$s_\lambda(-y\,/\,x)$ in \cite{macd}.) It follows that the expansion
of $s_\lambda(x\,/\,y)$ in terms of power sums is given by
  $$ s_\lambda(x\,/\,y) = \sum_{\nu\vdash n}z_\lambda^{-1}
  \chi^\lambda(\nu) \left( p_\nu(x)-p_\nu(y)\right). $$
Hence from equation~(\ref{eq:flasum}) we obtain the following result.

  \begin{theorem} \label{thm:vkla}
\emph{Let $\lambda^n\vk (a;b)$. Then for a fixed partition $\alpha$ we
  have} 
  $$ f^{\lambda^n/\alpha} = f^{\lambda^n}s_\alpha(a\,/-\! b)
  (1+O(1/n)). $$ 
  \end{theorem}

An explicit statement of Theorem~\ref{thm:vkla} does not seem to have
been published before. However, it was known by Vershik and Kerov and
appears in the unpublished doctoral 
thesis of Kerov. It is also a simple consequence of Okounkov's
formula \cite[Thm.\ 8.1]{ok} for $f^{\lambda/\alpha}$ in terms of
shifted Schur functions. The asymptotics of shifted Schur functions is
carried out (in slightly greater generality) in \cite[Thm.\ 8.1 and
Cor.\ 8.1]{k-o-o}. A special case of Theorem~\ref{thm:vkla} appears in
\cite[Thm.\ 1.3]{r-v}.

Theorem~\ref{thm:vkla} can be made more explicit in certain cases for
which the super-Schur function $s_\alpha(a\,/\,\!-\!b)$ can be
explicitly evaluated. In particular, 
suppose that $\alpha$ consists of an $i\times j$ rectangle with a
shape $\mu=(\mu_1,\dots,\mu_i)$ attached at the right and 
the conjugate $\nu'$ of a shape $\nu=(\nu_1,\dots,\nu_j)$ attached at
the bottom. Thus
  $$ \alpha=(\mu_1+j,\dots,\mu_i+j,\nu'_1,\nu'_2,\dots). $$
Then (e.g., \cite[pp.\ 115--118]{lit}\cite[(4) on p.\ 59]{macd})
  $$ s_\alpha(a_1,\dots,a_i\,/\,-b_1,\dots,-b_j) =s_\mu(a_1,\dots,a_i) 
  s_\nu(b_1,\dots,b_j)\prod_{i,j}(a_i+b_j). $$
In certain cases we can explicitly evaluate $s_\mu(a_1,\dots,a_i)$ or
$s_\nu(b_1,\dots, b_j)$, e.g., when $a_1=\cdots=a_i$ or
$b_1=\cdots=b_j$. See \cite[Thm.\ 7.21.2 and Exer.\ 7.32]{ec2}. 
Note also that when $\mu=\nu=\emptyset$ (so $\alpha=(j^i)$) we have
simply  
  $$  s_{(j^i)}(a_1,\dots,a_i\,/\,-b_1,\dots,-b_j)
  =\prod_{i,j}(a_i+b_j). $$ 

\textsc{Acknowledgement.} I am grateful to Andrei Okounkov for
providing much of the information about Theorem~\ref{thm:vkla}
mentioned in the paragraph following the statement of this theorem.


\end{document}